\documentclass[10pt]{article}
\usepackage{amsmath,amssymb}
\usepackage{graphicx,epsfig,epstopdf}
\usepackage{psfrag}
\usepackage{color}
\usepackage{amsfonts}
\usepackage{latexsym}
\usepackage{pstricks}
\usepackage{bbm}

\makeatletter\@addtoreset{equation}{section}\makeatother
\newtheorem{theorem}{\bf Theorem}
\newtheorem{remark}[theorem]{Remark}

\newcounter{constantsnumber}

\newcommand{\CQ}{\mathcal{Q}}

\newcommand{\CT}{\mathcal{T}}

\newcommand{\fb}{\mbox{\boldmath{$f$}}}

\newcommand{\bu}{\mbox{\boldmath{$u$}}}

\newcommand{\bv}{\mbox{\boldmath{$v$}}}
\newcommand{\bV}{\mbox{\boldmath{$V$}}}

\newcommand{\bx}{\mbox{\boldmath{$x$}}}

\newcommand{\bzero}{\mbox{$\bf 0$}}
\newcommand{\BBR}{\mbox{$\mathbb{R}$}}
\newfont{\twelvemsb}{msbm10 at 11.6pt}

\renewcommand{\div}{\mathop{\rm div}\nolimits}

\newcommand{\Span}{\mathop{\rm span}}

\usepackage{colordvi}

\title{A quadrilateral 'mini' finite element for the Stokes problem
 using a single bubble function}
\author{Bishnu P. Lamichhane
\thanks{School of Mathematical and Physical Sciences, University of Newcastle, Callaghan, NSW 2308,  {\tt Bishnu.Lamichhane@anu.edu.au}}}

\begin{document}
\maketitle

\begin{abstract}
We consider a quadrilateral 'mini' finite element for approximating the solution of 
Stokes equations using a quadrilateral mesh.
We use the standard bilinear  finite element space enriched with element-wise defined bubble functions 
  for the  velocity and the standard bilinear finite element space
for  the pressure space. With a simple modification of the standard bubble function
we show that a single bubble function is sufficient to ensure the inf-sup condition. 
We have thus improved an earlier result on the quadrilateral 'mini'  element, 
where more than one bubble function are used to get the stability.
\end{abstract}

\begin{keywords}
Stokes equations, mixed finite elements,  Mini finite element, inf-sup condition, bubble function
\end{keywords}

{\bf AMS subject classification}.
65N30, 65N15, 74B10

\pagestyle{myheadings}
\thispagestyle{plain}

\section{Introduction}
A very simple finite element method for the Stokes problem for a simplicial mesh  is 
presented by Arnold, Brezzi and Frotin \cite{ABF84}, where the velocity space is discretised by using 
the standard linear finite element space enriched with element-wise bubble functions and the pressure 
space is discretised by using the standard linear finite element space. 
The enrichment of the velocity space is done to ensure the stability of the finite element method, and 
this increases one vector degree of freedom 
per element. 
An extension of the finite element method to the quadrilateral mesh is done by Bai 
\cite{Bai97}, where the author enriches the velocity space with more than a single vector bubble function per element. 
The inf-sup condition is proved by using a macro element technique \cite{Ste90a}, where 
a single element is used as a macro element. 

In this article we show that with a small modification of the standard bubble function we can get 
the stability just by using a single vector bubble function per element. The main difference with 
the technique proposed by Bai \cite{Bai97} is that it is not possible to show the inf-sup condition using a single  element as a macro element. 
We need to use a macro element consisting of four elements to prove the inf-sup condition 
in our situation.  Another relevant finite element method is presented by Lamichhane \cite{Lam14}, where two different 
meshes are used to discretise the velocity and the pressure space, and a single vector bubble degree of freedom per element is used to get the stability. The pressure space is discretised by the space of piecewise 
constant functions on the dual mesh.  However, the main difficulty of the technique presented by 
Lamichhane \cite{Lam14} is that the bubble function is obtained by 
multiplying the standard bubble function by the gradient of a bilinear basis function, and hence the bubble function 
cannot be defined on  a reference element. The standard bubble function on the unit square 
 is  the lowest degree polynomial  which vanishes on the boundary of the square. 
 Here we modify the standard bubble function \cite{ABF84,Bai97} to get stability of the numerical scheme by using 
 a single vector bubble function per element  with a continuous pressure approximation. 
We also  investigate two choices of bubble functions, where both of them can be defined on a reference element. 
 Since the first mini finite element 
is introduced for simplicial meshes \cite{ABF84} with a single bubble function 
per element, this new contribution gives a unified 
framework for quadrilaterals and triangles. The idea can easily be extended to 
the three-dimensional case.

\section{Stokes equations}\label{sec:bvp}
This section is devoted to the introduction of the boundary 
value problem of the Stokes equations.  Let $\Omega$ in $\BBR^2$, be a 
bounded domain with polygonal  boundary $\Gamma$. 
For a prescribed body force $\fb \in [L^2(\Omega)]^2$, the Stokes 
equations with homogeneous Dirichlet boundary condition  in $\Gamma$ reads
\begin{equation}
\begin{array}{ccc}
-\nu\Delta \bu + \nabla p &=& \fb\quad\text{in}\quad \Omega \\
\div  \bu &=& 0 \quad\text{in}\quad \Omega 
\end{array}
\end{equation}
with $\bu=\bzero$ on $\Gamma$, 
where $\bu$ is the velocity, $p$ is the pressure, and $\nu$ denotes
the viscosity of the fluid. 

Here we use standard notations $L^2(\Omega)$, $H^1(\Omega)$ and 
$H^1_0(\Omega)$ for Sobolev spaces, see \cite{BS94,Cia78} for details. 
Let $\bV :=[H^1_0(\Omega)]^2$  be the vector Sobolev space with 
 inner product $(\cdot,\cdot)_{1,\Omega}$ and norm $\|\cdot\|_{1,\Omega}$
defined in the standard way:  $(\bu,\bv)_{1,\Omega} := \sum_{i=1}^2
(u_i,v_i)_{1,\Omega}$, and the norm being induced by this inner
product. We also define another subspace $M$ of $L^2(\Omega)$ as 
\[ P = \left\{ q \in L^2(\Omega): \, \int_{\Omega} q \, dx =0\right\}.\]

The weak formulation of the Stokes equations is to find $(\bu,p) \in   
\bV\times P$ such that 
\begin{equation} \label{stokesw}
\begin{array}{ccc}
\nu \int_{\Omega} \nabla\bu: \nabla \bv\,dx  &+
\int_{\Omega} \div v \,p\,dx &= \ell(\bv),\quad \bv \in \bV,\\
\int_{\Omega} \div \bu\, q \,dx &&= 0,\quad q \in P,
\end{array}
\end{equation}
where $ \ell(\bv) = \int_{\Omega} \fb\cdot \bv\,dx.$
It is well-known that the weak formulation of the Stokes problem is
well-posed \cite{GR86}. In fact, if the domain $\Omega$ is convex, and 
$\fb \in [L^2(\Omega)]^2$, we 
have $\bu \in [H^2(\Omega)]^2$, $p \in H^1(\Omega)$ and 
the a priori estimate holds 
\[ \|\bu \|_{2,\Omega}  + \|p \|_{1,\Omega}  \leq C \|\fb \|_{0,\Omega},\]
where the constant $C$ depends on the domain $\Omega$.

 \section{Finite element discretizations}
We consider a quasi-uniform triangulation $\CT_h$  of the 
polygonal domain $\Omega$, where $\CT_h$
consists of convex quadrilaterals.
The finite element meshes are defined by maps from the
reference square $\hat{K}=(0,1)^2$.

Let $\CQ_1(\hat K)$ be the space of bilinear  polynomials in $\hat K$. 
We start with the finite element space of 
continuous functions whose restrictions to an element
$K$ are obtained by maps of {\em bilinear} functions from the reference element: 
\begin{equation}
S_h:= \left \{v_h \in H_{0}^1(\Omega),\
v_h|_K=\hat{v}_h\circ F^{-1}_{K}, \hat{v}_h\in
\CQ_{1}(\hat{K}),~~K\in \CT_h\right \}, 
\end{equation}
where $F_K: \hat K \to K$ is an iso-parametric map. 
We note that the iso-parametric map $F_K$ is generated by using 
the basis functions of $\CQ_1(\hat K)$.
It is clear that if $\hat{v}\in \CQ_{1}(\hat{K})$, then
$\hat{v}\circ F^{-1}_{K}$ is in general not a polynomial on the
quadrilateral $K$. 

In the following we assume that each element $K \in \CT_h$ is a parallelogram and 
the map $F_K$ is affine. 
Let $b_K$  be a bi-variate polynomial of $\bx \in \BBR^2$ with $b_K =0$ on $\partial K$ and $b_K(\bx_K) =1$, 
where $\bx_K\in \BBR^2$ is  the centroid of $K$. This is called
a bubble function  corresponding to the element $ K \in \CT_h$.
Defining the space of bubble functions 
\begin{eqnarray}\label{bubble}
B_h := \{ b_h \in C^0(\Omega): 
b_h|_{K} = c_K b_K,\; c_K \in \BBR,\;  K \in \CT_h\}, 
\end{eqnarray}
we introduce our finite element space for  velocity as
$\bV_h = [S_h \oplus  B_h]^2$.
The finite element space for the pressure is taken as the standard 
bilinear  finite element space 
\begin{equation}\label{press}
S^*_h:= \left \{v_h \in L_0^2(\Omega) \cap H^1(\Omega),\
v_h|_K=\hat{v}_h\circ F^{-1}_{K}, \hat{v}_h\in
\CQ_{1}(\hat{K}),~~K\in \CT_h\right \}\ .
\end{equation}

Then, the finite element approximation of \eqref{stokesw}
is defined as a solution to the following problem: 
find $(\bu_h,p_h) \in \bV_h \times S^*_h$ such that
\begin{equation}\label{stokesd}
\begin{array}{llllllll}
a(\bu_h,\bv_h)+b(\bv_h,p_h)&=&\ell(\bv_h),\quad &\bv_h &\in& \bV_h , \\
b(\bu_h,q_h) &=&0, \quad &q_h &\in& S^*_h.
\end{array}
\end{equation}

We need the following conditions to prove that there is a unique solution of 
the discrete problem \eqref{stokesd} and the discrete solution converges optimally 
to the continuous solution. 
\begin{enumerate}
\item The bilinear forms 
$ a(\cdot,\cdot)$ on $\bV_h \times \bV_h$ and  $ b(\cdot,\cdot)$ on  $\bV_h \times S^*_h$ are continuous. 
 \item The bilinear form $ a(\cdot,\cdot)$ on $\bV_h \times \bV_h$ 
 is elliptic. 
 \item 
There exists a constant $\beta> 0$ independent 
of the mesh-size such that  for any $q_h \in S^*_h$, we have 
\begin{eqnarray}\label{infsup} 
\sup_{\bv_h \in \bV_h}
\frac{b(\bv_h,q_h)}{\|\bv_h\|_{1, \Omega}} \geq \beta \|q_h\|_{0,\Omega}.
\end{eqnarray}
The smallest constant $\beta$ with the property 
\begin{eqnarray}\label{infsupn} 
  \beta  = \inf_{q_h \in S_h^*} \sup_{\bv_h \in \bV_h}
\frac{b(\bv_h,q_h)}{\|\bv_h\|_{1, \Omega}, \|q_h\|_{0,\Omega}} 
\end{eqnarray}
is called the inf-sup constant.
\end{enumerate}
\section{The Macro-Element Technique} 
We  prove the inf-sup condition \eqref{infsup} using a 
macro-element technique proposed by Stenberg \cite{Ste90a}.
A macro-element $M$ is a connected set of elements in $\CT_h$. 
Moreover, two macro-elements $M_1$ and $M_2$ are said to be equivalent 
if they can be mapped continuously onto each other \cite{Ste90a}. 
We define the following three spaces associated with  the macro-element $M_i$:
\[ \bV_h^i =  [H^1_0(M_i)]^2 \cap \bV, \;
S^i_h= \left\{v_h \in H^1(M_i),\ v_h|_K=\hat{v}_h\circ F^{-1}_{K}, \hat{v}_h\in
\CQ_{1}(\hat{K}),~~K \in \CT_h,\ K\subset M_i \right\}, \] 
and \[ 
B_i = \left\{ q_h\in S^i_h|\, b(\bv_h,q_h)=0,\;\bv_h \in \bV_h^i\right\}.\]
Moreover, we denote by $\Gamma_h$ the set of all edges in $\CT_h$ interior to $\Omega$. 
The macro-element partition ${\mathcal{M}_h}$ of $\Omega$ then consists of 
macro-elements $\{M_i\}_{i=1}^N$ with $\bar\Omega = \bigcup_{i=1}^N \bar M_i$. 
The macro-element technique is given by the following theorem \cite{Ste90a}.
\begin{theorem}
Suppose that there is a fixed set of equivalence classes ${\mathcal{E}_j}$, $j=1,\cdots,q$, of macro-elements, a positive integer $L$, and a macro-element partition ${\mathcal{M}_h}$ such that 
\begin{enumerate}
\item[(M1)] For each $M_i \in {\mathcal{E}_j}$, $j=1,\cdots,q$, the space $B_i$ is one-dimensional, consisting of functions 
that are constant on $M_i$. 
\item[(M2)]  Each $M_i \in {\mathcal{M}_h}$ belongs to one of the classes  ${\mathcal{E}_j}$, $j=1,\cdots,q$. 
\item[(M3)]  Each $T \in \CT_h$ is contained in at least one and not more than $L$ macro-elements of  ${\mathcal{M}_h}$.
\item[(M4)]  Each $ e \in \Gamma_h$ is contained in the interior of at least one and not more than $L$ 
macro-elments of ${\mathcal{M}_h}$.
\end{enumerate}
Then the inf-sup condition \eqref{infsup} is satisfied. 
\end{theorem}
\begin{figure}[htb]
  \centerline{ \epsfig{figure=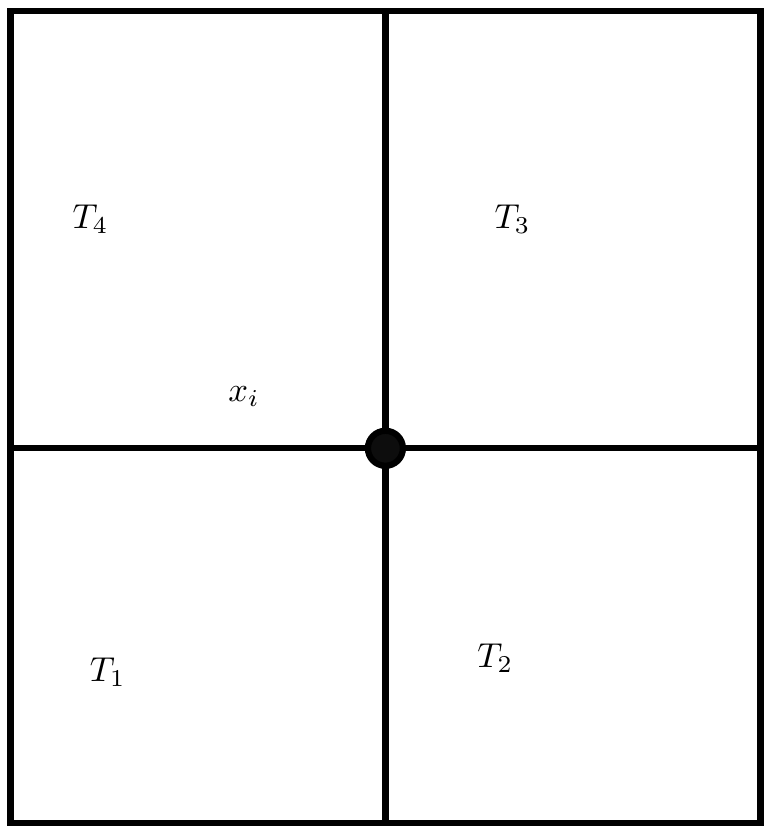,width=0.6\textwidth,height=0.4\textheight}}
 \caption{The set $M_i$, where four elements of $\CT_h$ touch the vertex $\bx_i$}
 \label{dualv}
\end{figure}

In the following we consider a macro-element consisting of four squares as shown in Figure 
\ref{dualv}. With this partition of macro-elements we can see that Assumptions (M2)--(M4) are all satisfied. 
We now show that the proof of Assumption (M1) depends on the choice of 
bubble functions. 
\subsection{Choice of bubble functions}
\begin{figure}[htb]
  \centerline{ \epsfig{figure=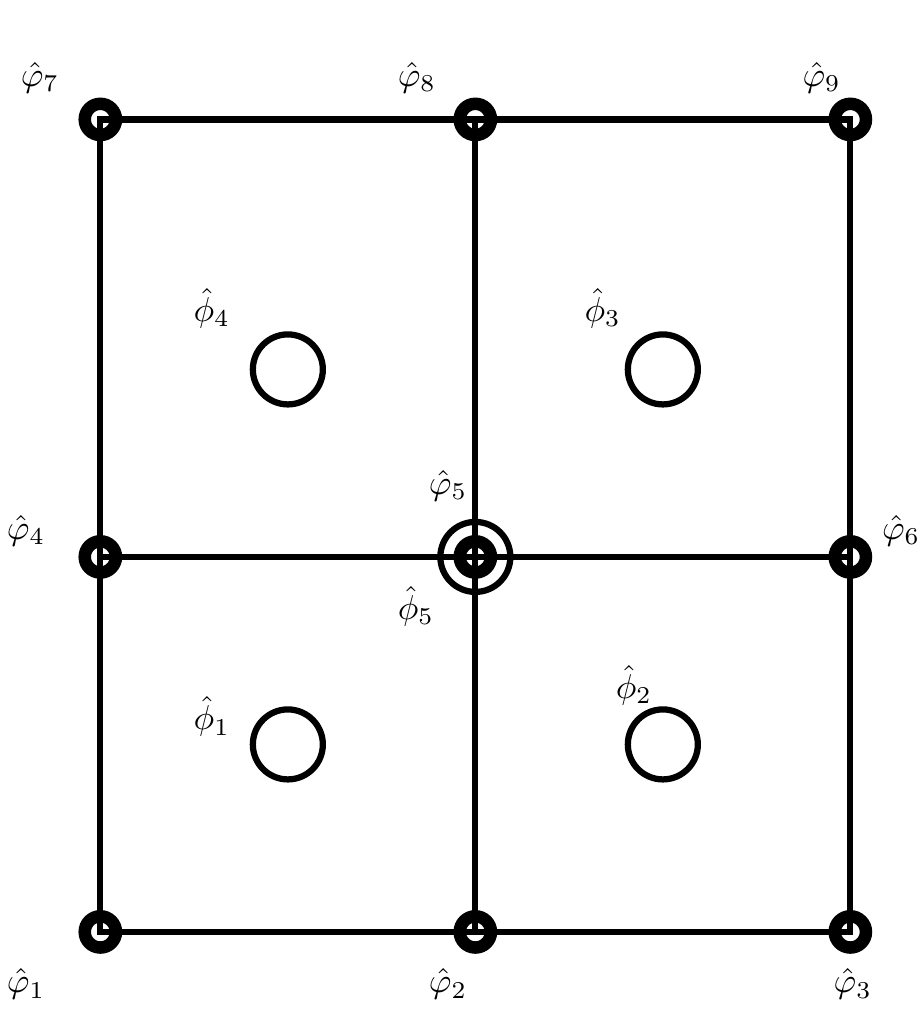,width=0.6\textwidth,height=0.4\textheight}}
 \caption{The numbering of functions $\hat\phi_k$ and $\hat \varphi_j$ on the reference square $\hat S$}
 \label{dualvn}
\end{figure}

For simplicity of calculation we assume that $M_i$ is a parallelogram so that there is an invertible affine mapping 
$F_i: \hat S \to M_i$, which transforms the square $\hat S=[-1,1]^2$ to $M_i$ with the property 
\begin{equation}\label{affine}
\begin{bmatrix}
x\\
y
\end{bmatrix}  =  A_i  \begin{bmatrix}
\xi\\
\eta
\end{bmatrix}  +  \begin{bmatrix}
x_0\\
y_0
\end{bmatrix},
\end{equation}
where $A_i$ is a 2 by 2 matrix,  $(x,y) \in M_i$ and $(\xi,\eta) \in \hat S$. 
Let $V_h^i=\Span  \{\phi_k\}_{k=1}^5$,  $ \bV_h^i = [V_h^i]^2$ 
and $S_h^i = \Span  \{\varphi_k\}_{k=1}^9$. We use the notation 
$\hat \phi_k$ and $\hat \varphi_k$ to denote corresponding basis functions 
on the square $\hat S$, where $\hat \phi_k$ and $\hat \varphi_k$  are functions 
of $\xi$ and $\eta$. We have shown the numbering of functions  $\hat\phi_k$ and $\hat \varphi_j$  
on the reference square $\hat S$ in Figure \ref{dualvn}, where 
we have used big circles for the functions in  $\bV_h$, and small 
circles for functions in $S_h^*$.

Let $\bv_h \in \bV_h^i$ with $\bv_h = \sum_{k=1}^5  \bv_k\phi_k $ and 
$\bv_k   \in \BBR^2$. Then 
\[ b(\bv_h, q_h) = 
\int_{M_i} \nabla \cdot \bv_h \, q_h \,dx  = 
\sum_{k=1}^5 \int_{M_i} \bv_k \cdot \nabla \phi_k \, q_h  \, dx.\]
Using a chain rule we write 
\[ \nabla \phi_k = A_i^{-T} \left(\hat\nabla \hat \phi_k \circ F_i^{-1}\right),\]
where $\hat \nabla$ denotes the gradient on the reference square $\hat S$. 
Let $q_h = \sum_{j=1}^9 q_j \varphi_j$, and thus 
\[ \int_{M_i} \nabla \cdot \bv_h \, q_h \,dx   = \sum_{k=1}^5 \sum_{j=1}^9q_j \bv_k \cdot 
\int_{M_i} \nabla \phi_k \, q_j\varphi_j  \, dx = |\det A_i| \sum_{k=1}^5 \sum_{j=1}^9q_j \bv_k \cdot 
\int_{\hat S} A_i^{-T} \hat\nabla \hat \phi_k \hat \varphi_j \, d\hat x.
\]
We see that we can find a matrix $\tilde D$ such that 
\[  \int_{M_i} \nabla \cdot \bv_h \, q_h \,dx  = \vec{q}^T \tilde D \vec{v},\]
where 
\[ \vec{q} = \begin{bmatrix}
q_1\\
q_2\\
\vdots\\
q_9
\end{bmatrix},\quad\text{and} \quad 
\vec{v} = \begin{bmatrix}
\bv_{1}\\
\bv_{2}\\
\vdots\\
\bv_{5}
\end{bmatrix}=  \begin{bmatrix}
v_{1}\\
v_{2}\\
\vdots\\
v_{10}
\end{bmatrix}.\]
Thus we need to show that the rank of the matrix $\tilde D$ is 8 in order to prove that 
the dimension of the space $B_i$ is one. 

Since $A_i$ is an invertible matrix, the rank of the matrix will be unchanged if we replace 
$M_i$ by the reference element $\hat S$, so that we want to investigate the 
rank of the matrix $D\in \BBR^{10 \times 9}$, where the $j$th row of  $D$ is 
\[ \left[\int_{\hat S} \partial_{\xi} \hat\phi_1\hat\varphi_j\,d\hat x, 
\int_{\hat S} \partial_{\eta} \hat\phi_1\hat\varphi_j\,d\hat x,
\int_{\hat S} \partial_{\xi} \hat\phi_2\hat\varphi_j\,d\hat x,
\int_{\hat S} \partial_{\eta} \hat\phi_2\hat\varphi_j\,d\hat x, \cdots, 
\int_{\hat S} \partial_{\xi} \hat\phi_5\hat\varphi_j\,d\hat x,
\int_{\hat S} \partial_{\eta} \hat\phi_5\hat\varphi_j\,d\hat x\right].\]

\subsubsection{Standard bubble functions}
Consider the unit square $ K = (0,1)^2$ in two dimensions. 
We start with the standard choice of 
the bubble function $b_K=16 xy(1-x)(1-y)$. 
The matrix $D$ is explicitly computed as 
\[ D= 
 \left[ \begin {array}{cccccccccc} \frac{2}{9}&\frac{2}{9}&0&0&0&0&0&0&\frac{1}{12}&\frac{1}{12}
\\ \noalign{\medskip}-\frac{2}{9}&\frac{2}{9}&\frac{2}{9}&\frac{2}{9}&0&0&0&0&0&\frac{1}{3}
\\ \noalign{\medskip}0&0&-\frac{2}{9}&\frac{2}{9}&0&0&0&0&-\frac{1}{12}&\frac{1}{12}
\\ \noalign{\medskip}\frac{2}{9}&-\frac{2}{9}&0&0&0&0&\frac{2}{9}&\frac{2}{9}&\frac{1}{3}&0
\\ \noalign{\medskip}-\frac{2}{9}&-\frac{2}{9}&\frac{2}{9}&-\frac{2}{9}&\frac{2}{9}&\frac{2}{9}&-\frac{2}{9}&\frac{2}{9}&0&0
\\ \noalign{\medskip}0&0&-\frac{2}{9}&-\frac{2}{9}&-\frac{2}{9}&\frac{2}{9}&0&0&-\frac{1}{3}&0
\\ \noalign{\medskip}0&0&0&0&0&0&\frac{2}{9}&-\frac{2}{9}&\frac{1}{12}&-\frac{1}{12}
\\ \noalign{\medskip}0&0&0&0&\frac{2}{9}&-\frac{2}{9}&-\frac{2}{9}&-\frac{2}{9}&0&-\frac{1}{3}
\\ \noalign{\medskip}0&0&0&0&-\frac{2}{9}&-\frac{2}{9}&0&0&-\frac{1}{12}&-\frac{1}{12}\end {array}
 \right].
\]
We compute the rank of this matrix using {\sc maple} and obtain it to be $7$. 
Thus in this case the dimension of the space $B_i$ will be two. Hence there is no hope of 
getting the inf-sup condition for this choice of the bubble function.

\subsubsection{The first choice of bubble functions}
In the next step, we consider the bubble function 
\[ b_K = 64 \varphi_K  xy(1-x)(1-y),\]
where $\varphi_K $ is the standard bilinear basis function 
corresponding to the lower-left corner of the square $K$. 
Since $\varphi_K = (1-x)(10y)$, 
the bubble function $b_K$ on the reference square $K$ can be defined as 
\[ b_K = 64 (1 - x)(1 - y)xy(1-x)(1-y). \]
Defined in this way the bubble function $b_K$ does not depend on the local 
numbering of the vertices of $K$.  In this case, the matrix $D$ has rank 8, and is computed as 
\[ D = \left[ \begin {array}{cccccccccc} {\frac {4}{15}}&{\frac {4}{15}}&0&0
&0&0&0&0&\frac{1}{12}&\frac{1}{12}\\ \noalign{\medskip}-{\frac {4}{15}}&{\frac {8}{45}
}&{\frac {4}{15}}&{\frac {4}{15}}&0&0&0&0&0&\frac{1}{3}\\ \noalign{\medskip}0&0
&-{\frac {4}{15}}&{\frac {8}{45}}&0&0&0&0&-\frac{1}{12}&\frac{1}{12}
\\ \noalign{\medskip}{\frac {8}{45}}&-{\frac {4}{15}}&0&0&0&0&{\frac {
4}{15}}&{\frac {4}{15}}&\frac{1}{3}&0\\ \noalign{\medskip}-{\frac {8}{45}}&-{
\frac {8}{45}}&{\frac {8}{45}}&-{\frac {4}{15}}&{\frac {4}{15}}&{
\frac {4}{15}}&-{\frac {4}{15}}&{\frac {8}{45}}&0&0
\\ \noalign{\medskip}0&0&-{\frac {8}{45}}&-{\frac {8}{45}}&-{\frac {4}
{15}}&{\frac {8}{45}}&0&0&-\frac{1}{3}&0\\ \noalign{\medskip}0&0&0&0&0&0&{
\frac {8}{45}}&-{\frac {4}{15}}&\frac{1}{12}&-\frac{1}{12}\\ \noalign{\medskip}0&0&0&0
&{\frac {8}{45}}&-{\frac {4}{15}}&-{\frac {8}{45}}&-{\frac {8}{45}}&0&
-\frac{1}{3}\\ \noalign{\medskip}0&0&0&0&-{\frac {8}{45}}&-{\frac {8}{45}}&0&0
&-\frac{1}{12}&-\frac{1}{12}\end {array} \right].
\]

\begin{remark}
We have used the gradient of the bilinear function $\varphi_K$ 
to construct  a vector bubble function associated with the element $K$ in 
\cite{Lam14}. Since  the construction  of the bubble function 
using the gradient of $\varphi_K$ cannot be done on a reference element, 
this new bubble function is computationally much easier.
\end{remark}

\subsubsection{The second choice of bubble functions}
It is interesting to see if we can multiply the bubble function by a linear function   and 
obtain the stability. For this purpose we can choose a bubble function on the unit square 
$(0,1)^2$ as 
\[ b_K = (a+bx+cy) xy(1-x)(1-y),\quad abc \neq 0.\]
For simplicity we choose 
\[ b_K= 8 (1+x+y) xy(1-x)(1-y).\]
We note that the factor 8 is used to force the value of the bubble function at the centroid of the square to be 1. 
The resulting matrix $D$ has also rank 8 in this case, and hence the 
dimension of the space $B_i$ is one. Moreover, the matrix $D$ is computed as 
\[ D =  \left[ \begin {array}{cccccccccc} {\frac {19}{90}}&{\frac {19}{90}}&0
&0&0&0&0&0&\frac{1}{12}&\frac{1}{12}\\ \noalign{\medskip}-{\frac {19}{90}}&{\frac {7}{
30}}&{\frac {19}{90}}&{\frac {19}{90}}&0&0&0&0&0&\frac{1}{3}
\\ \noalign{\medskip}0&0&-{\frac {19}{90}}&{\frac {7}{30}}&0&0&0&0&-\frac{1}{12}&\frac{1}{12}\\ \noalign{\medskip}{\frac {7}{30}}&-{\frac {19}{90}}&0&0&0&0&
{\frac {19}{90}}&{\frac {19}{90}}&\frac{1}{3}&0\\ \noalign{\medskip}-{\frac {7
}{30}}&-{\frac {7}{30}}&{\frac {7}{30}}&-{\frac {19}{90}}&{\frac {19}{
90}}&{\frac {19}{90}}&-{\frac {19}{90}}&{\frac {7}{30}}&0&0
\\ \noalign{\medskip}0&0&-{\frac {7}{30}}&-{\frac {7}{30}}&-{\frac {19
}{90}}&{\frac {7}{30}}&0&0&-\frac{1}{3}&0\\ \noalign{\medskip}0&0&0&0&0&0&{
\frac {7}{30}}&-{\frac {19}{90}}&\frac{1}{12}&-\frac{1}{12}\\ \noalign{\medskip}0&0&0&0
&{\frac {7}{30}}&-{\frac {19}{90}}&-{\frac {7}{30}}&-{\frac {7}{30}}&0
&-\frac{1}{3}\\ \noalign{\medskip}0&0&0&0&-{\frac {7}{30}}&-{\frac {7}{30}}&0&0
&-\frac{1}{12}&-\frac{1}{12}\end {array} \right] .\]

\begin{remark}
The proof of stability is presented for the two-dimensional case. However, 
this can be extended to the three-dimensional case without a major change. 
\end{remark}

\begin{remark}
It is interesting to see if we can use a quadratic  function symmetric about the centroid of the element 
to multiply the standard bubble function. 
To check this we use a bubble function on the reference square $K = (0,1)^2$ defined as 
\[ b_{K} = xy\left(x^2+y^2-x-y+\frac{33}{2}\right)(1-x)(1-y),\]
and compute the matrix $D$. In this case, the rank of the matrix $D$ is just 7, and hence the dimension of 
the space $B_i$ is 2. 
\end{remark}

An immediate consequence of the above discussion is the well-posedness 
of the discrete problem \eqref{stokesd}. From the theory of saddle point 
problem, see, e.g., \cite{BF91}, we have the following theorem.
\begin{theorem} \label{th1}
The discrete problem \eqref{stokesd} has exactly one solution 
$(\bu_h,p_h) \in \bV_h \times S^*_h$, which is uniformly stable
with respect to the data $\fb$, and there exists a 
constant $C$ independent of the mesh-size $h$ such that
\begin{eqnarray*}
\|\bu_h\|_{1,\Omega} + \|p_h\|_{0,\Omega} \leq C\|\fb\|_{0,\Omega}.
\end{eqnarray*}
\end{theorem}
The convergence theory is provided by an abstract result about the 
approximation of saddle point problems,
see \cite{BF91}.
\begin{theorem} \label{th2}
Assume that $(\bu,p)$ and $(\bu_h,p_h)$
be the solutions of problems \eqref{stokesw} and \eqref{stokesd},
respectively. Then, we have the following error estimate:
\begin{multline} \label{se}
\|\bu-\bu_h\|_{1,\Omega}+\|p-p_h\|_{0,\Omega}  \leq C 
\left(\inf_{\bv_h \in \bV_h}\|\bu-\bv_h\|_{1,\Omega}+
\inf_{q_h \in S^*_h}\|p-q_h\|_{0,\Omega}\right).
\end{multline}
\end{theorem}

\section{Numerical Results}\label{sec:num}

In this section we present two numerical experiments to verify the 
optimal a priori error estimate and some numerical experiments to verify the 
inf-sup condition for the proposed 
finite element scheme. 
For both examples we consider a 
simple unit square $\Omega = (0,1)^2$. 
\subsection{Verify a priori error estimate} 
For both examples we consider  a uniform 
initial triangulation consisting of four squares.

\paragraph{First example.}
 For the first example we choose the exact solution $\bu = (u_1,u_2)$ as 
 \[ u_1 = - 2\, x^2\, y\, \left(2\, y - 1\right)\, {\left(x - 1\right)}^2\, \left(y - 1\right),\quad
u_2 = 2\, x\, y^2\, \left(2\, x - 1\right)\, \left(x - 1\right)\, {\left(y - 1\right)}^2.\]
We use the kinematic viscosity $\nu =1$. The exact solution for the 
pressure is chosen as 
\[  p = x(1-x)(1-2y),\]
so that $ p \in L^2_0(\Omega)$. 
The exact solution $\bu$ satisfies the homogeneous Dirichlet boundary condition on $\partial \Omega$, and 
the right hand side function $\fb$ is computed by using the exact solution $\bu$ and the pressure $p$.
We have presented the errors in 
the velocity and the pressure approximation using the  $H^1$-norm and 
the $L^2$- norm, respectively in 
Table \ref{exa1b1} for the first choice of the bubble function, and 
in Table \ref{exa1b2} for the second choice of the bubble function. 
We note that the standard choice of the bubble function leads to a singular matrix.
From the presented tables we can see 
the optimal convergence of the 
velocity approximation in the $H^1$ and $L^2$-norms, and 
a super-convergence result for the pressure in the 
$L^2$-norm. As we expect a convergence rate of order 
$1$  for the pressure approximation in the $L^2$-norm but 
get a better approximation of order $1.5$, 
this is a super-convergence. This 
better convergence is due to the fact that 
 we have used the standard continuous bilinear finite element space 
 for the pressure approximation. 
 We can also observe that all errors are smaller for the second choice of bubble functions.

\begin{table}[!htb]
\caption{Discretization errors for the velocity and pressure, Example 1 (First choice)}
\begin{center}
\begin{tabular}{|c|c|c|c|c|c|c|c|c|}\hline
level $l$ & \# elem. & \multicolumn{2}{c|}{$\|u -u_h\|_{1,\Omega}$} & \multicolumn{2}{c|}{$\|u -u_h\|_{0,\Omega}$}
 & \multicolumn{2}{c|}{$\|p -p_h\|_{0,\Omega}$} \\ \hline
1 &    16 & 3.23129e-02 & & 3.03116e-03 & & 1.76150e-02 &  \\\hline
 2 &    64 & 1.58286e-02 & 1.03 & 8.24246e-04 & 1.88 & 7.00356e-03 & 1.33  \\\hline
 3 &   256 & 7.79938e-03 & 1.02 & 2.06421e-04 & 2.00 & 2.50753e-03 & 1.48  \\\hline
 4 &  1024 & 3.87699e-03 & 1.01 & 5.12144e-05 & 2.01 & 8.78516e-04 & 1.51  \\\hline
 5 &  4096 & 1.93346e-03 & 1.00 & 1.27289e-05 & 2.01 & 3.08875e-04 & 1.51  \\\hline
 6 & 16384 & 9.65545e-04 & 1.00 & 3.17131e-06 & 2.00 & 1.08856e-04 & 1.50  \\\hline

\end{tabular}
\end{center}
\label{exa1b1}
\end{table}

\begin{table}[!htb]
\caption{Discretization errors for the velocity and pressure, Example 1 (Second choice)}
\begin{center}
\begin{tabular}{|c|c|c|c|c|c|c|c|c|}\hline
level $l$ & \# elem. & \multicolumn{2}{c|}{$\|u -u_h\|_{1,\Omega}$} & \multicolumn{2}{c|}{$\|u -u_h\|_{0,\Omega}$}
 & \multicolumn{2}{c|}{$\|p -p_h\|_{0,\Omega}$} \\ \hline
 1 &    16 & 3.16876e-02 & & 2.89325e-03 & & 1.17765e-02 &  \\\hline
 2 &    64 & 1.56503e-02 & 1.02 & 7.90369e-04 & 1.87 & 4.31789e-03 & 1.45  \\\hline
 3 &   256 & 7.75922e-03 & 1.01 & 1.99983e-04 & 1.98 & 1.44890e-03 & 1.58  \\\hline
 4 &  1024 & 3.86716e-03 & 1.00 & 4.99365e-05 & 2.00 & 4.93948e-04 & 1.55  \\\hline
 5 &  4096 & 1.93102e-03 & 1.00 & 1.24544e-05 & 2.00 & 1.71287e-04 & 1.53  \\\hline
 6 & 16384 & 9.64934e-04 & 1.00 & 3.10849e-06 & 2.00 & 5.99594e-05 & 1.51  \\\hline

\end{tabular}
\end{center}
\label{exa1b2}
\end{table}

\paragraph{Second example.}
For the second example we consider an exact solution 
given in \cite{BDG06}, where the exact solution for the velocity 
$\bu = (u_1,u_2)$ is given by 
\[ u_1 = x+{x}^{2}-2\,xy+{x}^{3}-3\,x{y}^{2}+{x}^{2}y,\quad 
u_2 = -y-2\,xy+{y}^{2}-3\,{x}^{2}y+{y}^{3}-x{y}^{2},\]
and the exact solution for the pressure is given by 
\[ p = xy+x+y+{x}^{3}{y}^{2}-\frac{4}{3}.\]
We use the  kinematic viscosity $\nu =1$ and the exact solution 
to compute the right-hand side function $\fb$.
As in the first example we compute the errors in 
the velocity and the pressure approximation using the  $H^1$--norm and 
the $L^2$- norm, respectively. The numerical results are tabulated in 
Table \ref{exa2b1} and  \ref{exa2b2} for the two choices of 
bubble functions, respectively. 
As in the first example, we can see the optimal convergence rates for the 
velocity approximation in $H^1$ and $L^2$-norms, and a better 
convergence rate for the pressure in $L^2$-norm. 
We also observe that all errors are smaller for the second choice of 
bubble functions although the difference is quite small in this example.

\begin{table}[!htb]
\caption{Discretization errors for the velocity and pressure, Example 2 (First choice)}
\begin{center}
\begin{tabular}{|c|c|c|c|c|c|c|c|c|}\hline
level $l$ & \# elem. & \multicolumn{2}{c|}{$\|u -u_h\|_{1,\Omega}$} & \multicolumn{2}{c|}{$\|u -u_h\|_{0,\Omega}$}
 & \multicolumn{2}{c|}{$\|p -p_h\|_{0,\Omega}$} \\ \hline
 1 &    16 & 6.96126e-01 & & 3.33821e-02 & & 2.25132e+00 &  \\\hline
 2 &    64 & 3.39100e-01 & 1.04 & 8.37772e-03 & 1.99 & 5.58680e-01 & 2.01  \\\hline
 3 &   256 & 1.66684e-01 & 1.02 & 2.09556e-03 & 2.00 & 1.59539e-01 & 1.81  \\\hline
 4 &  1024 & 8.26546e-02 & 1.01 & 5.24458e-04 & 2.00 & 4.49273e-02 & 1.83  \\\hline
 5 &  4096 & 4.11633e-02 & 1.01 & 1.31193e-04 & 2.00 & 1.28191e-02 & 1.81  \\\hline
 6 & 16384 & 2.05425e-02 & 1.00 & 3.28081e-05 & 2.00 & 3.80370e-03 & 1.75  \\\hline
\end{tabular}
\end{center}
\label{exa2b1}
\end{table}

\begin{table}[!htb]
\caption{Discretization errors for the velocity and pressure, Example 2 (Second choice)}
\begin{center}
\begin{tabular}{|c|c|c|c|c|c|c|c|c|}\hline
level $l$ & \# elem. & \multicolumn{2}{c|}{$\|u -u_h\|_{1,\Omega}$} & \multicolumn{2}{c|}{$\|u -u_h\|_{0,\Omega}$}
 & \multicolumn{2}{c|}{$\|p -p_h\|_{0,\Omega}$} \\ \hline
 1 &    16 & 6.96024e-01 & & 3.23184e-02 & & 5.93926e+00 &  \\\hline
 2 &    64 & 3.35337e-01 & 1.05 & 7.82819e-03 & 2.05 & 4.04732e-01 & 3.88  \\\hline
 3 &   256 & 1.65795e-01 & 1.02 & 1.97572e-03 & 1.99 & 6.07983e-02 & 2.73  \\\hline
 4 &  1024 & 8.24467e-02 & 1.01 & 4.97135e-04 & 1.99 & 1.78268e-02 & 1.77  \\\hline
 5 &  4096 & 4.11137e-02 & 1.00 & 1.24714e-04 & 2.00 & 5.88206e-03 & 1.60  \\\hline
 6 & 16384 & 2.05304e-02 & 1.00 & 3.12328e-05 & 2.00 & 1.98964e-03 & 1.56  \\\hline
\end{tabular}
\end{center}
\label{exa2b2}
\end{table}

\section{Conclusion}
In this contribution we present a finite element method for Stokes equations using 
continuous bilinear finite elements enriched with 
bubble functions for the velocity approximation and continuous bilinear finite elements for 
the pressure. In contrast to an earlier contribution we show that 
a single vector bubble function per element is enough to 
guarantee the stability of the discrete linear system. 
The numerical results also demonstrate the optimal convergence rates 
for the velocity and pressure approximation. 

\section*{Acknowledgement}
Support from the near miss  grant of the University of Newcastle 
 is gratefully acknowledged.

\bibliographystyle{elsart-num-sort}
\bibliography{total}

\begin{thebibliography}{10}
\expandafter\ifx\csname url\endcsname\relax
  \def\url#1{\texttt{#1}}\fi
\expandafter\ifx\csname urlprefix\endcsname\relax\def\urlprefix{URL }\fi

\bibitem{ABF84}
D.~Arnold, F.~Brezzi, M.~Fortin, A stable finite element for the {S}tokes
  equations, Calcolo 21 (1984) 337--344.

\bibitem{Bai97}
W.~Bai, A quadrilateral 'mini' finite element for the {S}tokes problem,
  Computer Methods in Applied Mechanics and Engineering 143 (1997) 41--47.

\bibitem{BDG06}
P.~Bochev, C.~Dohrmann, M.~Gunzburger, Stabilization of low-order mixed finite
  elements for the {S}tokes equations, SIAM Journal on Numerical Analysis 44
  (2006) 82--101.

\bibitem{BS94}
S.~Brenner, L.~Scott, The Mathematical Theory of Finite Element Methods,
  Springer--Verlag, New York, 1994.

\bibitem{BF91}
F.~Brezzi, M.~Fortin, Mixed and hybrid finite element methods,
  Springer--Verlag, New York, 1991.

\bibitem{Cia78}
P.~Ciarlet, The Finite Element Method for Elliptic Problems, North Holland,
  Amsterdam, 1978.

\bibitem{GR86}
V.~Girault, P.-A. Raviart, Finite {E}lement {M}ethods for {N}avier-{S}tokes
  {E}quations, Springer-Verlag, Berlin, 1986.

\bibitem{Lam14}
B.~Lamichhane, A mixed finite element method for nearly incompressible
  elasticity and stokes equations using primal and dual meshes with
  quadrilateral and hexahedral grids, Journal of Computational and Applied
  Mathematics 260 (2014) 356--363.

\bibitem{Mal81}
D.~Malkus, Eigenproblems associated with the discrete {LBB} condition for
  incompressible finite elements, International Journal of Engineering Science
  19 (1981) 1299--1310.

\bibitem{Ste90a}
R.~Stenberg, Error analysis of some finite element methods for the stokes
  problem, Mathematics of Computation 54 (1990) 495--508.

\end{thebibliography}
\end{document}